\documentclass{amsart}

\usepackage[foot]{amsaddr}
\usepackage{graphicx}
\usepackage{hyperref}
\usepackage{natbib}
\bibpunct{(}{)}{;}{a}{}{,}
\renewcommand{\cite}{\citep}

\title{Mathematical Monsters}
\author[Andrew Aberdein]{Andrew Aberdein$^{*}$}
\address{$^{*}$School of Arts \& Communication, Florida Institute of Technology, Melbourne FL.}
\thanks{{Published in}: %Aberdein, Andrew (2019).
Compagna, Diego and Steinhart, Stefanie, edd., 
{\em Monsters, Monstrosities, and the Monstrous in Culture and Society}, pp.~391--412. Vernon Press, Wilmington, DE, 2019.}
\begin{document}
\maketitle

Monsters lurk within mathematical as well as literary haunts. 
I propose to trace some pathways between these two monstrous habitats. 
I start from Jeffrey Jerome Cohen's influential account of monster culture and explore how well 
mathematical monsters fit each of his seven theses \cite{Cohen96}.
The mathematical monsters I discuss are drawn primarily from three distinct but overlapping domains. I will describe these in much greater detail as they arise below, but here is a brief preview.
Firstly, late nineteenth-century mathematicians made numerous unsettling discoveries that threatened their understanding of their own discipline and challenged their intuitions. 
The great French mathematician Henri Poincar\'e characterised these anomalies as `monsters', a name that stuck. 
Secondly, the twentieth-century philosopher Imre Lakatos composed a seminal work on the nature of mathematical proof, in which monsters play a conspicuous role \cite{Lakatos76}. He reconstructs the emergence during the nineteenth century of a proof of the Euler Conjecture, which ascribes a certain property to polyhedra.%, 
\footnote{Polyhedra are three-dimensional figures with planar surfaces. Euler conjectured that $V - E + F = 2$, where $V$ is the number of vertices, $E$ the number of edges, and $F$ the number of faces. See \S\ref{sec:escapes} for further discussion.}
Lakatos coined such terms as `monster-barring' and `monster-adjusting' to describe strategies for dealing with entities whose properties seem to falsify the conjecture.
Thirdly, and most recently, mathematicians dubbed the largest of the so-called sporadic groups `the Monster', because of its vast size and uncanny properties, and because its existence was suspected long before it could be confirmed.% 
\footnote{Groups are mathematical structures of wide-ranging application; they may be most easily approached in terms of symmetry: different groups capture the symmetries of different sorts of object. Simple groups are groups that cannot be reduced to smaller groups---hence some authors call them `atoms of symmetry' \cite[8]{Ronan06}.
Sporadic groups are a special sort of finite simple group---ones that do not fit into the system of classification that accommodates most such groups.
See \S\ref{sec:desire} for further discussion.}

\section{`The monster's body is a cultural body'}\label{sec:body}
For Cohen, the monster is `an embodiment of a certain cultural moment---of a time, a feeling, and a place' \cite[4]{Cohen96}.
This poses two immediate questions: what, if anything, constitutes a cultural moment for mathematicians and can mathematics be embodied?
The philosopher Brendan Larvor argues that mathematical practices should be understood as cultural since `they are reproduced non-biologically, they are shared activities mediated by material artefacts, and they express norms and values' \cite[18]{Larvor16b}.
However, the question of \emph{what} constitutes a cultural moment for mathematicians is a fraught one. Is the culture in question the culture at large, as bolder historians and sociologists of mathematics would insist? 
Is it rather the culture of mathematicians' institutions, universities, learned societies, funding bodies and so forth, as other, more cautious sociologists and historians maintain?
Or is there a purely mathematical culture comprised of the mathematicians and their interactions, proceeding more or less independently of the wider world? 
Or perhaps many such cultures, divided by mathematical specialism or geographical contingency?
Or is the mathematician of no culture, his idiosyncrasies attributable only to his peculiar genius, to the psychology of the individual?

Some of the conflict between these various positions can be illustrated by a classic debate between a sociologist, David Bloor, and a philosopher, John Worrall, over the status of the mathematical monsters discussed by Lakatos.
Bloor, in his splendidly titled article `Polyhedra and the abominations of Leviticus' \cite{Bloor78} links the responses of nineteenth-century mathematicians to `monstrous' polyhedra to the fourfold classification of cultures devised by the anthropologist Mary Douglas: hierarchical, enclavist, individualistic, or isolate. As she explains,
\begin{quotation}
one is based on hierarchical community, and so in favour of formality and compartmentalization; 
one is based on equality within a group, and so in favour of spontaneity, and free negotiation, and very hostile to other ways of life;
one is the competitive culture of individualism;
and fourth is the culture of the isolate who prefers to avoid the oppressive controls of the other forms of social life \cite[42]{Douglas96}.%
\footnote{Bloor employs Douglas's earlier account of this taxonomy in terms of high or low `grid' and `group', where grid measures the strength of a society's internal hierarchies and group measures the strength of its barrier with the rest of the world. Hence hierarchists are high grid, high group; enclavists are low grid, high group; individualists are low grid, low group; and isolates are high grid, low group.}
\end{quotation}
Bloor traces the following sequence during the era in which proofs of Euler's Conjecture were under development:
`Eighteenth-century corporate universities' (enclavist) $\longrightarrow$
`Intended result of [Prussian university] reforms of 1812' (hierarchical) $\longrightarrow$
`Actual result of reforms by 1840: a competitive structure' (individualistic) \cite[265]{Bloor78}.
Each of these historical contexts corresponds to one of Douglas's four kinds of culture and is thereby associated with a method for dealing with monsters. Happily enough, the work which Lakatos saw as first exemplifying the most successful methodological response to the monster, which Bloor associates with an individualistic culture, appeared in 1847 (although in Munich, not Prussia, as Bloor concedes).

The difficulty with such a picture, of course, is that individual mathematicians do not necessarily fit. Mathematicians may share identical circumstances and yet go about their work in very different ways. 
Replying to Bloor, Worrall notes Poincar\'e's comment on his older compatriots Charles Hermite and Joseph Bertrand:
`They were scholars of the same school at the same time; they had the same education, they were under the same influences; and yet what a difference!' \citetext{\citealp[76]{Worrall79}, quoting \citealp[211]{Poincare13}}.
Of Poincar\'e himself, Worrall complains that `it just seems grossly implausible that a genius like Poincar\'e, who was engaged in critical and immensely fruitful debates on a variety of subjects with a variety of scientists, should have been simply reflecting his social circumstances in responding to polyhedra' \cite[78]{Worrall79}.
But there is nothing simple about it. Poincar\'e's genius led him to results that perhaps no other mathematician of his era could have found. Nonetheless, his social circumstances, including debates with his peers, helped to shape and direct that genius.
Likewise, Hermite and Bertrand may have had much in common, but they also had significant differences, not just in their temperaments but also in their respective networks of collaboration.% 
\footnote{For example, Hermite (but not Bertrand) took a close interest in contemporary German mathematics \cite[225]{Ferreiros16a}.}
Douglas's later work contains an implicit response to Worrall's challenge: the fourfold classification describes not only cultures, but also how individuals may accept or react against the prevailing culture \cite[44 ff.]{Douglas96}. 
Different mathematicians may be expected to accept or reject the various cultures of which they are a part to differing degrees. An applied mathematician, whose research is driven by the needs of science or industry, may reflect some cultural aspects more explicitly than a pure mathematician whose research is curiosity driven. A mathematician with many collaborators or students will belong to a culture of research in a different way from a more isolated contemporary.
But a mathematician who stood entirely outside of all such cultures would be a monster indeed.%
\footnote{By a further irony, Poincar\'e was himself once described as a `monster of mathematics' \cite[52]{Dieudonne73}. We might more naturally refer to him as a prodigy although, of course, `prodigy' also means monster.}

Worrall also suggests that `Bloor was perhaps misled by Lakatos's terminology: ``monster'' is, of course, a technical term in biology' \cite[78]{Worrall79}.
Indeed the appropriation of this biological terminology did not begin with Lakatos, but with Poincar\'e:
\begin{quotation}
Logic sometimes makes monsters. Since half a century we have seen arise a crowd of bizarre functions which seem to try to resemble as little as possible the honest functions which serve some purpose. No longer continuity, or perhaps continuity, but no derivative, etc. \dots\ 
If logic were the sole guide of the teacher, it would be necessary to begin with the most general functions, that is to say with the most bizarre. It is the beginner that would have to be set grappling with this teratologic museum (mus\'ee t\'eratologique) \cite[435 f.]{Poincare13}.%
\end{quotation}
But what is a teratologic museum?
Since the late eighteenth century, anatomists began to assemble collections of `monsters', developmentally anomalous human and animal foetuses. Poincar\'e's father was a professor of medicine, so we can assume that he had firsthand experience of such collections. More strikingly, 
`teratology' itself, an unambiguously technical term in biology that literally means `the study of monsters', was coined by Poincar\'e's grandfather-in-law.%
\footnote{The naturalist Isidore Geoffroy Saint-Hilaire introduced the term to describe the study of congenital malformations, a research field initiated by his father, Etienne Geoffroy Saint-Hilaire \cite[26]{Blumberg09}.
The two men were respectively grandfather and great-grandfather of Louise Poulain d'Andecy, Poincar\'e's wife \cite[20]{Rollet12}.}
So Lakatos and Poincar\'e are clearly using `monster' in the biological sense of an anomalous case that resists classification. That would be a problem for Bloor's application of Douglas (let alone my application of Cohen) if this sense of monster was unrelated to the latter's work. But, this is precisely what Douglas (and, as we shall see, Cohen) stresses the most \cite[see, for example,][126]{Douglas96}.

We have seen that mathematical practices take place within cultural contexts, but not that they are embodied. Surely mathematics is purely a matter of the mind?
In a recent paper entitled `The materiality of mathematics', the sociologist of mathematics Christian Greiffenhagen argues otherwise:
`writing (on paper, blackboards, napkins, beermats, or even in the air) and the development of representational techniques are indispensable for doing and thinking mathematics' \cite[505]{Greiffenhagen14}.
The truths of mathematics may be culture-independent, but the cultural process whereby mathematicians arrive at them and demonstrate their truth surely is not.
This distinction may be especially important for monsters. As objects of mathematics they have the same status as any other equally well-established mathematical object. But what makes them \emph{monsters} is their intersection with the specific mathematical culture that gave rise to them: `anomalies are not installed in nature but emerge from particular features of classificatory schemes' \cite[126]{Douglas96}.
Characteristically, they are monsters because they represent an affront to the intuitions of the mathematicians who discover them. Those intuitions are unavoidably cultural. In this sense, the monster's body is indeed a cultural body.

\section{`The monster always escapes'}\label{sec:escapes}
`No monster tastes of death but once', as Cohen has it \cite[5]{Cohen96}.
This is one of the most immediately recognisable traits of the monsters of film and folklore.
As Buffy says to Dracula `You think I don't watch your movies? You always come back' \cite[38$^{\prime}$00$^{\prime\prime}$]{Noxon00}.
Some mathematical monsters exhibit similar behaviour. For example, many of the monsters discussed by Lakatos were forgotten, rediscovered, forgotten again, re-rediscovered, reinterpreted, explained away, before they were finally put to good use. To see how this could happen we need to investigate the story Lakatos tells.

Polyhedra have been studied since antiquity.%
\footnote{Whereas polygons are plane figures bounded by straight edges, such as triangles, rectangles, or pentagons, polyhedra are their three-dimensional analogues: three-dimensional shapes with planar surfaces bounded by straight edges.
Regular polyhedra are polyhedra all of whose faces are regular (that is, equiangular and equilateral) polygons. They are five in number: the tetrahedron, comprised of four equilateral triangles; the cube, of six squares; the octahedron, of eight equilateral triangles; the dodecahedron, of twelve pentagons; and the icosahedron, of twenty equilateral triangles. 
Players of tabletop roleplaying games will recall these as the d4, d6, d8, d12, and d20 dice, respectively.
The proof that there can be only five such objects is pleasingly straightforward (ask yourself what the angles at the corners of each of the regular polygons are and how many of each could meet at a point). It is the last proposition of the last book of Euclid's \emph{Elements}, the most celebrated of mathematical textbooks \cite[XIII.18, 1190]{Heath06}.
Of course, not all polyhedra are regular. Roleplaying enthusiasts will also recall the d10, whose sides are kite-shaped, and perhaps the d30, with rhombus-shaped sides. More generally still, there are prisms, not all of whose sides are the same, and other figures assembled from various shapes in all sorts of ways.}
In the seventeenth century the mathematician and philosopher Ren\'e Descartes noted a property that all polyhedra, regular and irregular, seem to share. A century later the mathematician Leonhard Euler made the same observation, and put it into the form by which it is now best known:
$V - E + F = 2$, where $V$ is the number of vertices, $E$ the number of edges, and $F$ the number of faces.
The claim that this equation holds for every polyhedron is the Euler Conjecture. It clearly holds for tetrahedra ($4-6+4=2$), cubes ($8-12+6=2$), and the other regular polyhedra, and for many other examples besides. But just listing positive cases does not comprise a proof. A mathematical proof should establish that a theorem is true without exception. But as proofs of the conjecture emerged, so did apparent counterexamples. Throughout the nineteenth century numerous proofs of the conjecture were attempted, only to fall to counterexamples not anticipated by their authors. A final proof was not established until the very end of the century.

The Euler Conjecture may seem to be a mere mathematical curiosity. Nonetheless, the history of its attempted proofs is the subject of a masterpiece of the philosophy of mathematics:
\emph{Proofs and Refutations} is a singular book. It is structured as a dialogue, apparently taking place in a classroom, but in which the pupils, identified only by Greek letters, act as mouthpieces for the contrasting methodological tendencies of many of the nineteenth-century mathematicians who tangled with the Euler Conjecture. Meanwhile, the actual history unfolds in the footnotes.\footnote{Or not: even in the footnotes, Lakatos plays somewhat fast and loose with the historical record \cite[152]{Corfield03}.} As Lakatos reconstructs the story, the interruptions from monster polyhedra, apparent counterexamples to the conjecture, eventually led mathematicians to revise the concepts in which the conjecture was couched. Essentially, they came to realise that they had not really understood what a polyhedron was when they began searching for a proof. En route to this realisation, many different strategies are tried on for size.
These strategies worked to exclude, reinterpret, or otherwise kill off monster polyhedra. But the monsters kept escaping---or new monsters appeared to take their place.

\begin{figure}[htbp]
\begin{center}
\includegraphics[height=2cm]{%ArgMathExtra/
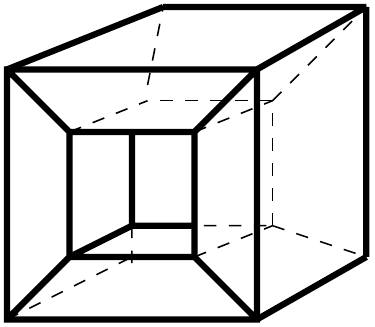}
\caption{The picture frame, fugitive non-Eulerian polyhedron?}
\label{fig:PictureFrame}
\end{center}
\end{figure}

The first monster to make an appearance is the hollow cube: `a solid bounded by a pair of cubes' \cite[13]{Lakatos76}. 
Since for each cube $V - E + F = 2$, then for both taken together as a single object $V - E + F = 4$.
In the historical record, this was first published in 1813 as a counterexample to Augustin Cauchy's 1811 `proof', then overlooked, then rediscovered in 1832.
In Lakatos's dialogue, the hollow cube provokes the following response:
\begin{quotation}
But why accept the counterexample? \dots\ Why should the theorem give way, when it has been proved? It is the `criticism' that should retreat. This pair of nested cubes is not a polyhedron at all. It is {a \emph{monster}, a pathological case, not a counterexample} \cite[14]{Lakatos76}.
\end{quotation}
This response, `monster-barring' to which I will return in \S\ref{sec:difference}, provokes
a series of monsters, for each of which a further monster-barring redefinition is proposed.
Many of these examples have a history of appearance, disappearance, and reappearance. For example, the picture frame ($16-32-16=0$), Fig.~\ref{fig:PictureFrame}, 
was first published alongside the hollow cube in 1813 (with a note from the journal's editor stating that he already knew of it), but was then overlooked in works published in 1827 and 1858 \cite[19]{Lakatos76}.
Other monster polyhedra included twin polyhedra (conjoined at a vertex or an edge), stellated polyhedra (with star-shaped faces), and the crested cube (a cube with a smaller cube stuck to one face).

Nor is monster-barring the only way to make monsters disappear. Lakatos also discusses monster-adjustment: where the monster-barrer revises definitions until the monster is no longer a polyhedron, the monster-adjuster revises definitions until the monster is no longer a monster. For the monster-adjuster,
\begin{quotation}
Monsters don't exist, only monstrous interpretations. One has to bar one's mind from perverted illusions, one has to learn how to see and how to define correctly what one sees. [Monster-adjustment] is therapeutic: where you---erroneously---`see' a counterexample, I teach you how to recognise---correctly---an example. I adjust your monstrous vision \cite[31]{Lakatos76}.
\end{quotation}
In some cases, the monster is merely an illusion, a misreading, in which case monster-adjustment may well be the best remedy. But the danger of this approach is that it does not explain, but only explains away. Where monsters are symptoms of deeper problems with the conjecture, monster-adjustment postpones the eventual recognition and remedy of these problems.

A corollary that Cohen draws from the fugitive nature of monsters is that `monstrous interpretation \dots\ must content itself with fragments' \cite[6]{Cohen96}.
This too is familiar from the culture of mathematical monsters. Consider these reflections on the Monster group (to which we shall return in \S\ref{sec:desire}):
\begin{quotation}
These things are so beautiful. It's such a pity that people can't see them. I mean, it's a kind of beauty that exists in the abstract, but we poor mortals will never see it. We can just get vague glimmerings. \dots\ 
Well, I mean, I know all the theorems. But there's still something that to me is unknown, unknowable. \dots\ Especially with the Monster, and I keep saying that it makes me sad that I'll probably never understand it \cite[264]{Roberts15a}.
\end{quotation}
This is John Conway, a mathematician who is as much of an expert on the Monster group as anyone.
His words may seem like false modesty, until one appreciates that the Monster group has more than $8 \times 10^{53}$ members. Where Conway has `vague glimmerings',
Cohen refers to `footprints, bones, talismans, teeth, shadows, obscured glimpses' as `signifiers of monstrous passing that stand in for the monstrous body itself' \cite[6]{Cohen96}. 
Mathematicians deal with different sorts of signifier from monster hunters, but they often have a similar task: to infer the properties of something they can only grasp in the most indirect fashion.
 
\section{`The monster is the harbinger of category crisis'}
Cohen's third thesis grounds his second: `The monster always escapes because it refuses easy categorization' \cite[6]{Cohen96}. This makes it `dangerous, a form suspended between forms that threatens to smash distinctions' \citetext{ibid.}. 
We have already seen some small-scale examples of this tendency at work in mathematics: the monster polyhedra disrupted early attempts to impose an easy categorization. Allusions to monsters also arise in much more fundamental areas of mathematical work.
An early example occurs in the origins of the calculus. Although Newton and Leibniz pioneered the study of rates of change of functions in the seventeenth century, their work relied upon the use of infinitesimal quantities, for which no consistent basis was then known. This was subject to criticism from an early date, but a consistent foundation was not presented until the nineteenth century.
The most astute early critic of the calculus, the philosopher and churchman George Berkeley, characterized infinitesimals in monstrous terms: they were the `ghosts of departed quantities' \cite[81]{Berkeley34}.
Much more recently, the infinitesimal has been rehabilitated, using modern mathematical techniques far beyond those available in the seventeenth century. Although the soundness of this work is beyond dispute, not all mathematicians find it worthwhile. Indeed, some of them have criticized the new infinitesimals in terms strikingly similar to those once used against the old---notably the distinguished French mathematician Alain Connes, who has referred to them as `chimeras' \cite[287]{Kanovei13}.
Although ghosts and chimeras are monsters of rather different pedigree, Berkeley and Connes are complaining about a property they share: resistance to categorization. The ghost is sometimes there and sometimes not; the chimera is a jumble of properties from different creatures.

A more fundamental category crisis arose in the nineteenth-century mathematics required for the rigorization of the calculus. Many concepts that had been hitherto given only informal definitions were replaced by formal explications, but it turned out---alarmingly---that the newly redefined concepts did not fit together in the same way as their pre-formal antecedents.
The irony was that dispatching one generation of monsters had given rise to another. As the mathematician and philosopher Solomon Feferman notes,
`one service that the monsters lurking around the corners provide is forcing us to don such armor [of rigorous proof] for our own protection. But if the proofs themselves produce such monsters, then the significance of \emph{what} is proved requires closer attention' \cite[328]{Feferman00}.

Crucially, mathematicians had assumed that continuity and differentiability were properties of more or less the same functions. Informally, a function is continuous if it doesn't have any `jumps'---that is, if its output only makes big changes when its input makes big changes---and it is differentiable whenever it has a derivative, that is when its rate of change can be determined. Functions that are differentiable at every point must be continuous. Intuitively, it may seem as though the converse should also be true, that continuous functions must be everywhere differentiable. This cannot be so, because there are some functions with isolated non-differentiable points, known as singularities, such as cusps or `sharp corners', that lack tangents. But it had been assumed that such functions would only have comparatively few, isolated singularities. Indeed, the French mathematician Andr\'e-Marie Amp\`ere published what he took to be a proof of this result in 1806.

\begin{figure}[htbp]
\begin{center}
\includegraphics[height=6cm]{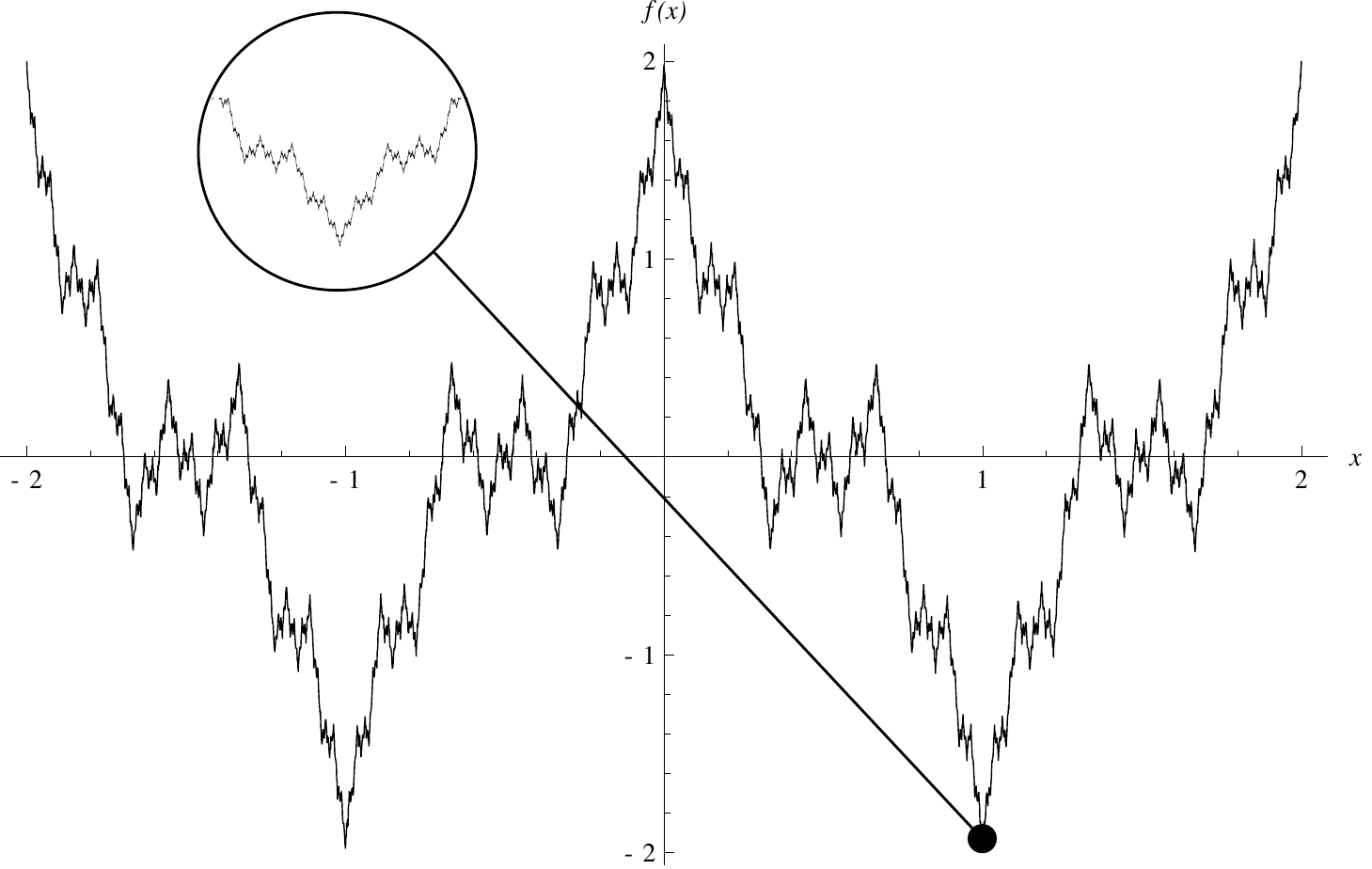}
\caption{Plot of Weierstrass function over the interval $[-2, 2]$, with detail enlarged. (Adapted from \url{https://en.wikipedia.org/wiki/Weierstrass_function}.)}
%\caption{Weierstrass's function $W$ with $a=.5$ and $b=5$ on $[0,3]$ \cite[22]{Thim03}. }
\label{fig:Weierstrass}
\end{center}
\end{figure}

However, starting in the 1830s, mathematicians began to turn up functions that were continuous but \emph{nowhere} differentiable \cite[for details, see][]{Thim03,Volkert08,Smorynski12}.
The Czech mathematician Bernard Bolzano was the first to demonstrate the existence of such a function, although his result remained in manuscript until long after his death (and long after similar results had been derived by others). The first published continuous but nowhere differentiable function is due to Karl Weierstrass, in a lecture given in 1872 and published in 1875. His function has the form
$$W(x)=\sum_{k=0}^{\infty}a^{k}\cos(c^{k}\pi x)$$
where $0<a<1$, $ab>1+3\pi/2$ and $b>1$ an odd integer. As may be seen from Fig.~\ref{fig:Weierstrass}, it is loosely sinusoidal, but between each summit the curve follows a smaller scale sinusoidal path, and between the summits of that path an even smaller scale sinusoidal path, and so on, ad infinitum. It is this last clause that is crucial: the curve stays just as detailed however much we magnify it. 

Self-similarity at ever smaller scales is a hallmark of an important class of fractals, a type of function which received a lot of attention in the late twentieth century, in part because increasing computer speeds made visualization of such functions feasible. Indeed the very word `fractal' is new---it was coined by Beno\^it Mandelbrot as recently as 1975.
While the classification of fractals is still an on-going project, and controversial in some details, they now have a clearly established place in mathematics. This was not the case at the time of Bolzano's and Weierstrass's discoveries, as Poincar\'e's animadversions against the exhibits of the `teratologic museum' demonstrate.
This suggests a point of contrast with the monsters analysed by Cohen, for whom `the monster's very existence is a rebuke to boundary and enclosure', its geography `always a contested cultural space' \cite[7]{Cohen96}. 
Progress in mathematics often consists in mapping that geography, resolving these contests, and constructing new monster-proof enclosures.
(Of course, new monsters tend to arise to test the new boundaries.)

\section{`The monster dwells at the gates of difference'}\label{sec:difference}
`The exaggeration of cultural difference into monstrous aberration is familiar enough', as Cohen reminds us \cite[7]{Cohen96}.
The most obvious mathematical counterpart to this xenophobic tendency in monster reception is Lakatos's method of monster-barring, which we encountered in \S\ref{sec:escapes}. He characterizes it, somewhat pejoratively, as follows:
\begin{quotation}
Using this method one can eliminate any counterexample to the original conjecture by a sometimes deft but always \emph{ad hoc} redefinition of the polyhedron, of its defining terms, or of the defining terms of its defining terms. We should somehow treat counterexamples with more respect, and not stubbornly exorcise them by dubbing them monsters \cite[23]{Lakatos76}.
\end{quotation}
Monster-barring is not confined to mathematics---similar techniques for winning arguments by manipulating definitions to exclude troublesome counterexamples are familiar from other contexts, under other names---some of which also allude to the xenophobic possibilities of monster-barring.
For example, the `{no true Scotsman}' move, in which the monster is barred by splitting the contested concept into `true' cases and others \cite[47]{Flew75}. Thus the chauvinistic Scotsman may reconcile his belief that `No true Scotsman takes sugar with his porridge' to the existence of an \emph{apparent} sugared-porridge-eating Scotsman, by dismissing the latter as no true Scotsman.

Bloor also brings out the scapegoating aspects of monster-barring, linking it to the anomaly response of
\begin{quotation}
Pollution-conscious societies [that] \dots\ survive by the threat of expulsion, or suffer repeated schism. They are frequently subject to outside threat and consequently their whole system of classification is pervaded by the dichotomy between the good inside and the evil and perverted outside. They need to exercise and symbolize high group control \cite[253]{Bloor78}.
\end{quotation}
Worrall complains of this link that 
`Bloor's argument \dots\ seems to depend quite heavily on there being an emotional component to ``monster-barring'', that it should involve ``turning in disgust'' from alleged counter-examples. But of course this emotional component is by no means essential' \cite[78]{Worrall79}.
Not essential, but by no means absent either: Lakatos quotes Hermite reacting to the functions discussed in the last section, `I turn aside with a shudder of horror from this lamentable plague of functions which have no derivatives' \cite[quoted in][19]{Lakatos76}.
It is significant that Lakatos is quoting from private correspondence, published posthumously; perhaps Hermite would have been more circumspect if writing for publication.
This reflects the sociologist Erving Goffman's distinction between the `front' and the `back' of cultural practices: the former intended for public consumption; the latter for private cooperation \cite[114]{Goffman71}.
Some commentators have found this reflected in mathematical practice. Thus Reuben Hersh distinguishes `mathematics in ``finished'' form, as it is presented to the public in classrooms, textbooks, and journals' from `mathematics as it appears among working mathematicians, in informal settings, told to one another in an office behind closed doors' \cite[128]{Hersh91}.
A proper understanding of mathematical culture must attend to both of these contexts but the front, intended as it is for a wider audience, is always the more accessible.
Bearing this practical distinction in mind, we may reflect that it is more remarkable that there is any trace of an emotional component to an apparently affectless problem-solving technique, than that it is not seen in every case.

Resort to monster-barring represents one crude conclusion that may be drawn from the presence of monsters: the monster must be excluded by whatever means necessary. 
Just as crudely, we could take the monsters to show that the whole basis of arriving at distinctions is radically unsound:
`By revealing that difference is arbitrary and potentially free-floating, mutable rather than essential, the monster threatens to destroy \dots\ the very cultural apparatus through which individuality is constituted and allowed' \cite[12]{Cohen96}.
A spokesman for this perspective with respect to mathematics was the philosopher Hans Hahn, for whom monsters such as continuous but nowhere differentiable functions had profound methodological implications:
`Again and again we have found that, even in simple and elementary geometric questions, intuition is a wholly unreliable guide. It is impossible to permit so unreliable an aid to serve as the starting point or basis of a mathematical discipline' \cite[1974]{Hahn33}.
For Hahn, this mandates 
`the expulsion of intuition from mathematical reasoning' and `the complete formalization of mathematics' \cite[1970]{Hahn33}.
The project Hahn envisages, of `completely formalizing mathematics, of reducing it entirely to logic' \cite[1971]{Hahn33}, was indeed a goal of much early twentieth-century foundational work in the philosophy of mathematics. But it foundered on deep conceptual difficulties and never much intersected with the research priorities of mainstream mathematicians. Subsequent thought has tended to less drastic solutions than the banishment of intuition. 
Feferman, for example, draws a much more modest conclusion:
\begin{quotation}
What, then, is one to say about the geometrical and topological monsters that are supposed to demonstrate the unreliability of intuition? The answer is simply that these serve as counterexamples to intuitively expected results when certain notions are used as explications which serve various purposes well enough but which do not have all expected properties \cite[322]{Feferman00}.
\end{quotation}
That is, the monsters show us the limits of our intuition, not that our intuition should be ignored. On the contrary, intuition has an enduring place in mathematical methodology, as does the practice of pushing it to its limits:
`it is standard mathematical practice to seek best possible results of an expected kind, and one way to achieve such is to make weakest possible assumptions on the given data. In this respect the mathematical monsters serve simply to provide counter-examples to further possible improvements' \citetext{ibid.}. This limitative role for the mathematical monster leads us directly to Cohen's fifth thesis.

\section{`The monster polices the borders of the possible'}
Marking limits is a further cultural use of monsters that Cohen identifies: `From its position at the limits of knowing, the monster stands as a warning against exploration of its uncertain demesnes' \cite[12]{Cohen96}.
This role has several dimensions with mathematical counterparts. We have seen one already, in Feferman's reflections on the limits of intuition.
Another is that particularly hard problems can be framed in monstrous terms, as a warning to mathematicians that their efforts may be in vain. The sociologist and historian of mathematics Donald MacKenzie notes how such notoriety attached to the four-colour conjecture, the innocent-seeming claim that four colours suffice to colour any planar map so that no adjacent regions are the same colour. Quoting Wolfgang Haken, one of the architects of the eventual proof (itself a computer-assisted monster, too vast for human inspection), MacKenzie states that the `conjecture had gained a reputation as a ``man-eating problem'' [whose] victims abandoned ordinary, tractable mathematics to spend years or even decades looking for a proof or a refutation' \cite[9]{MacKenzie99}. 
Indeed, MacKenzie's history of the four-colour conjecture invokes a monster in its title, `Slaying the Kraken', echoing the celebratory lines of the graph theorist (and occasional poet) Bill Tutte:
`Wolfgang Haken/Smote the Kraken/One! Two! Three! Four!/Quoth he: ``The monster is no more''' \cite[39 f.]{MacKenzie99}. When the conjecture was still open, Tutte had observed of it, `It is dangerous to work close to The Problem' \cite[quoted in][9]{MacKenzie99}.
The danger is a real one: mathematicians who devote all their effort to one big problem risk having nothing to show for it, to the detriment of their careers.

Cohen also observes that `policing the boundaries of culture [is] usually in the service of some notion of group ``purity''' \cite[15]{Cohen96}.
Mathematics has several notions of purity. 
The philosophers Mic Detlefsen and Andy Arana remark that
```Purity'' in mathematics has generally been taken to signify a preferred relationship between the resources used to prove a theorem or solve a problem and the resources used or needed to understand or comprehend that theorem or problem' \cite[1]{Detlefsen11}.
This is what is often termed `purity of method', the sense that an optimal proof should use only the resources of the field in which the problem arose, and not detour into other areas of mathematics. Such detours can be effective ways of tackling tough problems, but mathematicians often consider it worth pursuing a pure proof, even after an otherwise satisfactory proof has been established \cite[279 f.]{Dawson06}.
Another sense of purity is exhibited in the distinction between pure and applied mathematics. Here the purity emphasizes a boundary between mathematics and other scientific and technological disciplines: pure mathematics addresses problems that arise within mathematics; applied mathematics brings mathematical methods to bear on problems from elsewhere.
This boundary has always been as much political as objective: it served the institutional purposes of increasingly professionalized nineteenth-century mathematics departments to stress the purity of their work \cite[218]{Ferreiros16a}. Moreover, innovative scientists and engineers have often found very practical applications for mathematical results originally considered the purest of the pure: the centrality of prime factorization to modern cryptography is a standard example.
Taken to extremes, this sense of purity could become a pernicious ivory-tower mentality that could be used to excuse toleration of monstrosities in the political sphere: what the historian Herbert Mehrtens terms `irresponsible purity' in his survey of mathematics in Nazi Germany \cite{Mehrtens94}.

The boundary-patrolling function of monsters is also reflected in another of the methods of monster accommodation that Lakatos documents.
Rather than change the definitions so that the monster is excluded, as in monster-barring, a more productive approach may be to determine where the conjecture holds. This project of constructing a boundary to enclose the safe cases rather than to exclude the monsters, Lakatos calls exception-barring, or rather `a whole continuum of exception-barring attitudes' \cite[24]{Lakatos76}. He distinguishes primitive exception-barring, which just compartmentalizes the proof and its exceptions without any attempt to reconcile the two, from exception-barring strategies that are informed by the proof at issue. 
For these exception-barrers the proof at least provides `some inspiration for stating the conditions which determine a safe domain', and at best `a very fine delineation of the prohibited area' based upon `a careful analysis of the proof' \citetext{ibid.}.
Such an exception-barrer maintains that
\begin{quotation}
no conjecture is generally valid, but only valid in a certain restricted domain that excludes the \emph{exceptions}. I am against dubbing these exceptions `monsters' or `pathological cases'. That would amount to the methodological decision not to consider these as interesting \emph{examples} in their own right, worthy of a separate investigation \citetext{ibid.}.
\end{quotation}
This is a significant methodological improvement over monster-barring, but it is not the final word. The lingering problem is that addressing \emph{some} exceptions does not ensure that \emph{all} exceptions have been addressed. As Poincar\'e observed of such a manoeuvre, 
`We have put a fence around the herd to protect it from the wolves but we do not know whether some wolves were not already within the fence' \cite[quoted in][1186]{Kline72}.

Bloor sees exception-barring as characteristic of Douglas's hierarchical cultures: 
`large, diverse but stable system[s] of institutions', with
\begin{quotation}
an extensive repertoire of methods for responding to anomaly and for their reclassification. They will be fitted in somewhere, or the classificatory scheme will be expanded. Complicated rites of atonement; promotions and demotions; special exceptions; distinctions, assimilations and legal fictions will abound. And pervading the use of all these expedients there will be a vague sense of overriding unity \cite[254]{Bloor78}.
\end{quotation}
This neatly intersects with Cohen's remark that 
`The monster of prohibition \dots\ validated a tight, hierarchical system of naturalized leadership and control where every man had a functional place' \cite[13 f.]{Cohen96}.

\section{`Fear of the monster is really a kind of desire'}\label{sec:desire}
Cohen tells us that `The monster also attracts' \cite[16]{Cohen96}. This is clearly true of some mathematical monsters, at least for some mathematicians: the warnings against monster problems we encountered in the last section would be unnecessary if no one were drawn to such problems.
This thesis may be less obvious for the mathematical monsters we have discussed most, the anomalies and pathological cases. Certainly, some mathematicians, with Hermite, `turn aside with a shudder of horror'. But for other mathematicians
`the linking of monstrosity with the forbidden makes the monster all the more appealing as a temporary egress from constraint. \dots\ We distrust and loathe the monster at the same time we envy its freedom, and perhaps its sublime despair' \cite[17]{Cohen96}.
The concept of the sublime may be key to understanding this combination of fear and desire.
In one of the most influential analyses of the sublime, Edmund Burke characterizes it as follows:
\begin{quotation}
Whatever is fitted in any sort to {excite the ideas of pain and danger}, that is to say, whatever is in any sort terrible, or is conversant about terrible objects, or operates in a manner analogous to terror, is a source of the {\em sublime}; that is, it is productive of the strongest emotion which the mind is capable of feeling. When danger or pain press too nearly, they are incapable of giving any delight, and are simply terrible; but at certain distances, and with certain modifications, they may be, and they are delightful \cite[13 f.]{Burke56}.
\end{quotation}
Many informal reflections by mathematicians on the highs and lows of their professional lives could be framed in terms of the sublime, and sometimes these allusions are made explicit.
For example, the mathematician John Baez describes
`moments of exaltation that come from suddenly glimpsing a terrifyingly grand vista: sometimes shrouded in mist, sometimes lit by a lightning-bolt of insight' \cite{Baez08}.
This directly echoes Burke's observation that `{Greatness of dimension} is a powerful cause of the sublime' \cite[51]{Burke56}. 
Burke presumably intended `dimension' in the sense of scale---great extent along a given dimension---not in the sense of higher-dimensions, since the fourth dimension was not theorized until a century after his death. The Monster group, however, has greatness of dimension in both senses.
Its discovery has been described as `one of the most spectacular and mysterious mathematical achievements of the past fifty years' \cite[334]{Simons05}.

To gain some sense of what the Monster group is, and why it has attracted so much attention, we will need to venture briefly into group theory, the mathematical field in which it originates.
Technically, a group is a set, say $G$, acted on by an operation, say $*$, for which the following axioms hold:
\begin{description}
\item[Closure] $g*h \in G$ for all $g,h \in G$.
\item[Associativity] $(g*h)*k=g*(h*k)$.
\item[Identity] There is an identity element $e \in G$ such that $g*e = e*g = g$ for all $g \in G$.
\item[Inverse] For all $g \in G$ there exists $g^{-1}\in G$ such that $g*g^{-1} = g^{-1}*g = e$.
\end{description}
To see how this works, take a letter, say P, and think about some of the transformations we could perform upon it. For example, we could turn it upside down, \raisebox{7pt}{\reflectbox{\rotatebox{180}{P}}}, or back to front, \reflectbox{P}, or both at once, \raisebox{7pt}{\rotatebox{180}{P}}.
Counting leaving the original P alone as a null transformation, this gives us four transformations which comprise a group (known as the Klein four group, or $K_{4}$). We may think of them in terms of the distinct symbols that result from applying each to P: 
(P,
\raisebox{7pt}{\reflectbox{\rotatebox{180}{P}}},
\reflectbox{P},
\raisebox{7pt}{\rotatebox{180}{P}}).
These comprise the elements of this group. Its operation is just performing one transformation after another. For example, if we turn the letter upside down and then upside down again, it ends up back where it began. So this transformation is its own inverse (as are the other members of $K_{4}$), since leaving alone is clearly the identity element. With a bit more thought we can see that the elements are associative and closed, since they represent every possible combination of vertical and horizontal reflection: either, both, or neither.

As far as P is concerned, none of these transformations is a \emph{symmetry}, since P is clearly transformed. But suppose we had applied the transformations instead to X or H. Then (at least for sans serif fonts) the letters would remain unchanged. Hence the Klein four group represents (some of) the symmetries of these letters.
In contrast, consider a different set of four transformations, those resulting from a quarter-turn, a half-turn, a three-quarter-turn, and a full turn: 
(\raisebox{7pt}{\rotatebox{270}{P}},
\raisebox{7pt}{\rotatebox{180}{P}},
{\rotatebox{90}{P}},
P).
This is a different four-element group, the cyclic group of order four, $Z_{4}$. We can see that it is different by observing that H would not be symmetrical under this group (two of the operations would give us {\rotatebox{90}{H}}). But we can also see that the two groups have two elements in common, (\raisebox{7pt}{\rotatebox{180}{P}}, P). This is also a group, the cyclic group of order two, $Z_{2}$. The existence of this subgroup means that neither $K_{4}$ nor $Z_{4}$ are \emph{simple}.
However, the only subgroups of $Z_{2}$ are trivial: $Z_{2}$ itself and the identity element on its own. Hence $Z_{2}$ is simple, as is $Z_{p}$ for any prime number $p$. Thus the cyclic groups of prime order comprise a family of infinitely many simple groups.

There are several other families of simple groups, also with infinitely many members. One of the most impressive results of late twentieth-century mathematics was the classification of finite simple groups, which comprises literally thousands of pages of work, written by scores of mathematicians over several decades. This classification shows that every simple group of finite order must belong to one of these families, or to a small cluster of outlaws: the sporadic groups \cite[334]{Simons05}. These are finite simple groups that do not belong to any of the families.
In the 1860s and 70s, the French mathematician \'Emile L\'eonard Mathieu discovered the five sporadic groups that bear his name (listed here with their orders):
 \[ \begin{array}{cl}
    M_{11} & 2^{4}\cdot3^{2}\cdot5\cdot11 = 7,920 \\ 
    M_{12} & 2^{6}\cdot3^{3}\cdot5\cdot11 = 95,040 \\ 
    M_{22} & 2^{7}\cdot3^{2}\cdot5\cdot7\cdot11 = 443,520 \\ 
    M_{23} & 2^{7}\cdot3^{2}\cdot5\cdot7\cdot11\cdot23 = 10,200,960\\ 
    M_{24} & 2^{10}\cdot3^{3}\cdot5\cdot7\cdot11\cdot23 = 244,823,040\\ 
  \end{array}\]
Until the 1930s it remained controversial whether the Mathieu groups existed and, if so, whether they were really simple \cite[131]{Ronan06}. That no further sporadic groups were found in this period may have increased suspicions. But in 1965 the first of a brief spate of sporadic groups turned up.
Several of the new sporadic groups were actually smaller than some of the Mathieu groups, but some of them were very big indeed, notably {the Baby Monster}, which has order
\[2^{41} \cdot3^{13} \cdot5^{6} \cdot7^{2} \cdot11 \cdot13 \cdot17\cdot19\cdot23\cdot31\cdot47 = \]\[4154781481226426191177580544000000\] 
\[\approx 4 \times 10^{33}\]
and
{the Monster}, of order
\[2^{46} \cdot3^{20} \cdot5^{9} \cdot7^{6} \cdot11^{2} \cdot13^{3} \cdot17\cdot19\cdot23\cdot29\cdot31\cdot41\cdot47\cdot59\cdot71 = \]\[808017424794512875886459904961710757005754368000000000\] 
\[\approx 8 \times 10^{53}\]

The Monster is monstrous not only for its size, but for its uncanny properties. 
Specifically, it turns out to have profound connections to other, seemingly quite unrelated areas of mathematics. The discovery of these connections began as an apparent coincidence: a close relationship between two series of large numbers, the `character degrees' of the Monster and the coefficients of the $j$-function, which expresses a property of complex numbers \cite[192]{Ronan06}.
This, and other unexpected connections between the Monster and number theory, were dubbed `monstrous moonshine' by Conway, one of the pioneers of this field, from the slang term `for ``insubstantial or unreal'', ``idle talk or speculation'', ``an illusive shadow'' \dots\ to give the impression that matters here are dimly lit, and that [this technique] is ``distilling information illegally''' \cite[7]{Gannon06a}.
With collaborator Simon Norton, Conway published a series of conjectures extrapolating from the idea that these similarities could not just be coincidence \cite{Conway79}.
The eventual proof of these conjectures by Richard Borcherds, a former student of Conway, earned Borcherds a Fields medal, the most prestigious award in mathematics \cite[225]{Ronan06}.

It is worth noting that, although $8 \times 10^{53}$ is by many measures a big number, it is not especially large by the standards of many mathematicians. Much larger finite numbers are employed in some fields, such as combinatorics, some of them requiring innovative notations just to write down, since they are practically inexpressible in the standard notation of exponentiation.
Even these numbers are negligible by comparison with infinity itself. As Burke observed, infinity
`has a tendency to fill the mind with that sort of {delightful horror}, which is the most genuine effect, and truest test of the {sublime}' \cite[52]{Burke56}. 
For much of mathematical history it was accepted wisdom that the infinite could only be potential, not actual. That is, infinity might represent a never-to-be-attained limit, but it could not actually be assigned to completed sets of things. That perspective was overturned in the late nineteenth century by the German mathematician Georg Cantor, who devised an account of sets of infinite size. This would not be all that interesting a theory if there was only one infinite number, but Cantor showed that this cannot be so. 

With finite numbers, we often neglect to distinguish ordinal numbers (such as first, second, third) from cardinal numbers (such as one, two, three). With infinite numbers the distinction becomes much more important. 
The natural numbers $\{1,2,3,\dots\}$ are indeed equinumerous with (have the same cardinal number as) many other sets of numbers including the integers $\{\dots,-3,-2,-1,0,1,2,3,\dots\}$, the rationals (which can all be expressed as ratios of natural numbers) and the algebraic numbers (which are all expressible as roots of polynomial equations). That is, they all have the same cardinal number, $\aleph_{0}$.
But the set of real numbers, that is every number that can be expressed as a decimal expansion (possibly itself infinite in length), has demonstrably more members than there are natural numbers. In other words, it is uncountable. Of course, once we have two infinite numbers, we can start looking for more. 
Since Cantor's day, ever more remote families of infinite cardinals have been discovered, so much so that set theorists now routinely distinguish `large cardinals' from comparative minnows such as those studied by Cantor. 
Many of these families have exotic or grandiose names, although the mathematician Rudy Rucker observes that
\begin{quotation}
It is interesting to note that the smaller large cardinals have much grander names than the really big ones. Down at the bottom you have the self-styled inaccessible and indescribable cardinals loudly celebrating their size, while above, one of the largest cardinals quietly remarks that it is measurable, and the largest cardinals known simply point out past themselves with the comment that they are extendible \cite[265]{Rucker82}.
\end{quotation}
This observation is somewhat undercut by further developments: beyond the extendible lie almost huge, huge, and superhuge cardinals \cite[331]{Kanamori09}.

The development of a mathematics of the infinite was not without controversy. Some mathematicians regarded such work as intrinsically fanciful: internally consistent, perhaps, but devoid of any meaningful application. Some even went as far as to propose that such work should be jettisoned altogether.
This provoked a famous riposte from the doyen of German mathematics,
David Hilbert: `No one shall drive us out of the paradise which Cantor has created for us' \cite[191]{Hilbert26}.
This controversy finds an echo in Cohen's remark that
`The habitations of the monsters \dots\ are also realms of happy fantasy, horizons of liberation' \cite[18]{Cohen96}.

\section{`The monster stands at the threshold\dots\ of becoming'}
Cohen concludes his foray into monster culture by stating that `These monsters ask us how we perceive the world \dots\ They ask us why we have created them' \cite[20]{Cohen96}.
In \S\ref{sec:body} we encountered Poincar\'e invoking biological terminology to complain of the teratologic museum confronting the mathematician.
Modern biology has gone beyond teratology to \emph{teratogeny}, `the art of producing monstrous forms' \cite[254]{Blumberg09}. The creation of monsters---and attending to the questions that they have to ask---is also key to their role in mathematical discovery.

\emph{Proofs and Refutations} contains exactly one reference to biological literature \cite[22]{Lakatos76}.
Lakatos cites the early twentieth-century biologist Richard Goldschmidt, 
remembered for his now discredited idea of monsters as the source of new species, and particularly for his reference to 
`hopeful monsters, monsters which would start a new evolutionary line if fitting into some empty environmental niche' \cite[547]{Goldschmidt33}. 
Goldschmidt was reacting against a growing consensus in evolutionary biology that only external, environmental factors were relevant to the understanding of organisms. From that perspective, monsters were irrelevant, since they were presumed to be genetic dead-ends.
Much more recently evolutionary developmental biology (evo-devo, for short) has rehabilitated monsters---as sources of information, not of new species. As Pere Alberch, one of the originators of evo-devo puts it,
\begin{quotation}
I propose that it may be advantageous to turn to non-functional, grossly maladapted, teratologies when studying the properties of internal factors in evolution. These major deviations from normal development result in forms that are often lethal, and always significantly less well adapted than their progenitors. Therefore, one expects monsters to be consistently eliminated by selection. This is a useful property because if, in spite of very strong negative selection, teratologies are generated in a discrete and recurrent manner, this order has to be a reflection of the internal properties of the developmental system \cite[28]{Alberch89}.
\end{quotation}
Lakatos draws a similar moral:
`if we want to learn about anything really deep, we have to study it not in its ``normal'', regular, usual form, but in its critical state, in fever, in passion. 
If you want to know the normal healthy body, study it when it is abnormal, when it is ill.
If you want to know functions, study their singularities. 
If you want to know ordinary polyhedra, study their lunatic fringe' \cite[23]{Lakatos76}.
This insight lies behind the methodology that he judges most successful of all those he surveys. He calls this the method of proofs and refutations, but it has also been dubbed `monster-assimilating' \cite[108]{Caneva81}. 
The crucial insight is that monsters should be welcomed as a means of deepening proofs. A monster-assimilator reacts to a counterexample to a proof by adding an extra lemma to the proof that the counterexample refutes. This form of proof analysis leads to a revised conjecture and ultimately a deeper grasp of the concepts contained within the original conjecture \cite[127]{Lakatos76}.
For the Euler conjecture, the extra lemmas specify that there must be a way of travelling along edges between any pair of vertices (refuted by the hollow cube) and the surface must be simply connected (any loop we draw on it could be tightened to a point; refuted by the picture frame, Fig.~\ref{fig:PictureFrame}). But the proof analysis yields further dividends: we may generalize the result by dropping these lemmas; in doing so, $V-E+F$ will not always be $2$; instead it becomes an \emph{Euler characteristic} that classifies many different sorts of surface. Surfaces with Euler characteristic $2$ are those which are topologically equivalent to a sphere (that is, they could be `inflated' into spheres), whereas surfaces with Euler characteristic $0$, such as the picture frame, are topologically equivalent to a torus, or bagel.
Thus, by attending closely to monsters, a conjecture which might have seemed a mere mathematical curiosity motivates a profound and wide-ranging classification.

Bloor links such successful deployment of monster-assimilating to Douglas's individualistic cultures:
\begin{quotation}
We must ask what social forms exert a pressure towards innovation and novelty, and encourage transactions across the boundaries of existing classificatory schemes, dissolving them in change? Where is discontinuity more desired than regularity? Where can mistakes be tolerated and risks taken? Where is the tension most acutely felt between the missing of opportunities due to a reluctance to change, and missing them because of lack of sustained application? Which societies embody this contradiction in their very structure? The answer is: individualistic, pluralistic, competitive, and pragmatic social forms \cite[256]{Bloor78}.
\end{quotation}
These are the cultural conditions in which monsters may receive their due.

\end{document}